\documentclass[10pt]{amsart}
\usepackage{amsfonts}
\usepackage{graphicx}
\usepackage{psfrag}
\usepackage{palatino}
\usepackage{url}
\usepackage{amsmath,amsthm,amsfonts,amscd, amssymb, epsf} 
\usepackage{geometry}

\input psfig.sty
 \theoremstyle{plain}    
 \newtheorem{thm}{Theorem}[section]
 
 \numberwithin{equation}{section} 
 \theoremstyle{plain}
 \theoremstyle{plain}    
 \theoremstyle{definition}
 \newtheorem{defn}[thm]{Definition}
 \newtheorem{as}[thm]{Assumption} 
 \theoremstyle{plain}    
 \newtheorem{prop}[thm]{Proposition} 
 \newtheorem{rem}[thm]{Remark} 
 \newtheorem{lem}[thm]{Lemma}
 \newtheorem{cor}[thm]{Corollary}
 \newtheorem*{cor*}{Corollary}
 \newtheorem*{conj*}{Conjecture}

\newcommand{\bl}{\begin{lem}}
\newcommand{\el}{\end{lem}}

\newcommand{\bml}{\begin{multline}}
\newcommand{\eml}{\end{multline}}

\newcommand{\beq}{\begin{equation}}
\newcommand{\eeq}{\end{equation}}
\newcommand{\bp}{\begin{prop}}
\newcommand{\ep}{\end{prop}}
\newcommand{\bd}{\begin{defn}}
\newcommand{\ed}{\end{defn}}
\newcommand{\pf}{\begin{proof}}
\newcommand{\epf}{\end{proof}}

\newcommand{\tdv}{two-dimensional vector }
\newcommand{\td}{two-dimensional }
\newcommand{\thd}{three-dimensional  }
\newcommand{\fd}{finite-dimensional  }

\newcommand{\thds}{three-dimensional scalar }
\newcommand{\thdv}{three-dimensional vector }
\newcommand{\field}[1]{\ensuremath{\mathbb{#1}}}
\newcommand{\CC}{\field{C}}
\newcommand{\df}{\equiv}
\newcommand{\NN}{\field{N}}
\newcommand{\QQ}{\field{Q}}

\newcommand{\HHHH}{\field{H}}
\newcommand{\hh}{\field{H} \,}

\DeclareMathOperator*{\lto}{\mathnormal{o}(1)}

\DeclareMathOperator{\lsum}{\mathnormal{\sum}}
\DeclareMathOperator*{\psum}{ \lsum^\prime}
\DeclareMathOperator{\gip}{\Gamma ^ \prime _\infty   }
\DeclareMathOperator{\gi}{\Gamma_\infty   }
\DeclareMathOperator{\sz}{\mathnormal{Z}(\mathnormal{s},\Gamma,\chi)}
\DeclareMathOperator{\galp}{\Gamma ^ \prime _\alpha   }
\DeclareMathOperator{\gal}{\Gamma_\alpha   }
\DeclareMathOperator{\HH}{\HHHH^3}

\newcommand{\PP}{\field{P}}

\newcommand{\RR}{\field{R}}

\newcommand{\ZZ}{\field{Z}}

\newcommand{\F}{\mathcal{F}}
\newcommand{\E}{\mathcal{E}}
\newcommand{\scz}{\mathcal{S}}
\newcommand{\pg}{\mathcal{P}}
\newcommand{\K}{\mathcal{K}}
\newcommand{\D}{\mathcal{D}}
\newcommand{\hil}{\mathcal{H}}

\newcommand{\mC}{\mathcal{C}}

\newcommand{\smat}{\mathfrak{S}}

\DeclareMathOperator{\rep}{Rep(\Gamma,\mathnormal{V})}
\DeclareMathOperator{\PSL}{PSL}

\DeclareMathOperator{\SL}{SL}

\DeclareMathOperator{\vol}{vol}

\DeclareMathOperator{\R}{Re}
\DeclareMathOperator{\I}{Im}
\DeclareMathOperator{\pc}{PSL(2,\CC)}

\DeclareMathOperator{\ID}{id}
\DeclareMathOperator{\CUSP}{cusp}
\DeclareMathOperator{\PAR}{par}
\DeclareMathOperator{\CE}{ce}
\DeclareMathOperator{\NCE}{nce}
\DeclareMathOperator{\LOX}{lox}
\DeclareMathOperator{\cuspi}{\mathcal{CE}}

\DeclareMathOperator{\lds}{\mathnormal{ \frac{\phi^{\prime}}{\phi}}}

\DeclareMathOperator{\tr}{tr}

\DeclareMathOperator{\en}{\mathnormal{\mathcal{E}(T)}}
\DeclareMathOperator{\oen}{\mathnormal{\left|\mathcal{E}(T) \right|}}

\DeclareMathOperator{\ren}{\mathnormal{\mathcal{E}(R)}}

\newcommand{\cusps}{ \alpha \in \{ 1 \dots \kappa \}   }
\newcommand{\sing}{l \in \{ 1 \dots k_\alpha \}   }

\DeclareMathOperator{\cinf}{\mathnormal{\PP}}

\DeclareMathOperator{\hs}{\mathnormal{\hil(\Gamma,\chi)}}
\DeclareMathOperator{\lp}{\mathnormal{\Delta}}

\newcommand{\Z}{\mathbb{Z}}

\renewcommand{\Re}{\operatorname{Re}}

\newcommand{\ra}{\rightarrow}

\newcommand{\union}{\cup}

\begin{document}
\title[The Selberg trace formula and Selberg zeta-function]{The Selberg trace formula and Selberg zeta-function for  cofinite Kleinian groups with finite-dimensional unitary representations}
\author{Joshua S. Friedman}
\address{ Department of Mathematics, Stony Brook University, Stony Brook, NY 11794 }
\email{joshua@math.sunysb.edu}
\maketitle

\begin{abstract}
For cofinite Kleinian groups,  with finite-dimensional unitary representations, we derive the Selberg trace formula.  As an application we define the corresponding Selberg zeta-function and compute its divisor, thus generalizing results of Elstrodt, Grunewald and Mennicke  to non-trivial unitary representations.  We show that the presence of cuspidal elliptic elements sometimes adds ramification point to the zeta function.  In fact, if \mbox{$\mathcal{O} = \ZZ[-\frac{1}{2}+\frac{\sqrt{-3}}{2}]$} is the ring of Eisenstein integers, then the Selberg zeta-function of $\PSL(2,\mathcal{O})$ contains ramification points and is the sixth-root of a meromorphic function.
\end{abstract}
\tableofcontents{}
\section{Introduction}
The Selberg theory (Selberg trace formulas, Selberg zeta-functions, and related applications) has been well studied in both the two-dimensional scalar case (\cite{Iwaniec}) and the two-dimensional vector case\footnote{The works \cite{Roelcke}, \cite{Hejhal}, \cite{Fischer}, contain not only the two-dimensional vector case, but also its generalization, the case of unitary multiplier 
systems of arbitrary real weight. } (\cite{Roelcke}, \cite{Venkov}, \cite{Hejhal}, \cite{Fischer}).  By ``Two-dimensional vector case'' we mean: cofinite Fuchsian groups with finite-dimensional unitary representations, and the ``Scalar case'' refers to the case with the trivial  representation.  Elstrodt, Grunewald and Mennicke extended the Selberg theory to the three-dimensional scalar case in \cite{Elstrodt}.  By the ``Three-dimensional case'' we mean: cofinite Kleinian groups.  The main goal of this paper is to extend the Selberg theory to the three-dimensional vector case.

In this paper we derive the Selberg trace formula for cofinite Kleinian groups\footnote{A  Kleinian groups is referred to in some texts as a discrete group of isometries acting on hyperbolic three-space, or a discrete subgroup of $\pc.$},  with finite-dimensional unitary representations.  As an application we define the corresponding Selberg zeta-function and compute its divisor, thus generalizing results of Elstrodt, Grunewald and Mennicke \cite{Elstrodt} to non-trivial unitary representations. 

As one would expect, much of the \tdv and \thds cases extends in a straight forward manner to the \thdv case.  However, the extension of several parts of the Selberg theory are more subtle in the three-dimensional case, especially in the vector case.  One reason for this is because the set of finite-dimensional unitary representations of a fixed cofinite Kleinian groups is not well understood.  Another reason is related to the structure of the stabilizer subgroup of a cusp.  In the \td case the stabilizer subgroup of a cusp is a purely parabolic group that is isomorphic to a rank-one lattice,  while in three dimensions the stabilizer subgroup of a cusp is a non-abelian group that contains elliptic elements, with a finite-index purely parabolic subgroup that is  isomorphic to a rank-two lattice. The presence of elliptic elements in the stabilizer subgroup introduces some subtleties to the \thdv case, particularly in the computation of the divisor of the Selberg zeta-function.  In addition, the fact that the stabilizer subgroup (in the \thd case) contains a rank-two parabolic subgroup forced us to prove some additional estimates involving \td lattice sums.

A  Kleinian group  is a discrete subgroup of $\PSL(2,\CC) =   \SL(2,\CC) / \pm I. $ Each element of $\PSL(2,\CC) $ is identified with a M\"{o}bius transformation, and has a well-known action on hyperbolic three-space  $\HH$  and on its boundary at infinity$-$ the Riemann sphere $\PP^1$ (see \cite[Section 1.1]{Elstrodt})  . A Kleinian  group is \emph{cofinite} iff it has a fundamental domain $\F \subset \HH $ of finite hyperbolic volume. 

We use the following coordinate system for hyperbolic three-space, 
$\HH \df   \{(x,y,r)\in\RR^{3}~|~r>0 \} \df  \{ (z,r) ~| z \in \CC, ~r > 0 \} \df    \{z + rj\in\RR^{3}~|~r>0 \}, $ with the hyperbolic metric 
$$ ds^{2} \df \frac{dx^{2}+dy^{2}+dr^{2}}{r^{2}}, $$ and volume form 
$$ dv \df \frac{dx\, dy\, dz}{r^{3}}. $$
The Laplace-Beltrami operator is defined by
$$ \lp \df -r^{2}(\frac{\partial^{2}}{\partial x^{2}}+
\frac{\partial^{2}}{\partial y^{2}}+ \frac{\partial^{2}}{\partial r^{2}})+ r\frac{\partial}{\partial r}, $$ and it acts on the space of  smooth functions $f:\HH \mapsto V, $ where $V$ is a \fd complex vector space with inner-product $\langle~,~\rangle_V.$  

Suppose that $\Gamma $ is a cofinite Kleinian group and  $\chi \in \rep$ ($\rep$ is the space of \fd unitary representations of $\Gamma$ in $V$). Then the Hilbert space of \emph{$\chi-$automorphic} functions is defined by  
\begin{multline*}
\hs  \df  \{ f: \HH \ra V ~|~ f(\gamma P) = \chi(\gamma) f(P)~\forall \gamma \in \Gamma, \\ P \in \HH, $ and $\left<f,f \right> \df \int_{\F} \left<f(P),f(P)\right>_V\,dv(P) < \infty \}. 
\end{multline*}
Here $\F$ is a fundamental domain for $\Gamma $ in $\HH$,  and $\left<~,~\right>_V $ is the inner product on $V.$ Finally, let   $ \lp = \lp(\Gamma,\chi) $ be the corresponding positive self-adjoint Laplace-Beltrami operator on $\hs.$

Our first result is the spectral decomposition of $\lp$  on $\hs$ (see Theorem \ref{T:spectral}).  Except for one important point, the proof of the spectral decomposition theorem is analogous to the \tdv and \thds cases.  The one important point being, \emph{singularity at a cusp}.  To the best of the author's knowledge, prior to this paper, the notion of singularity at a cups was only defined for cofinite Fuchsian groups \cite{Selberg1} \cite{Venkov} \cite{Hejhal}.  In \S\ref{secSpectral} we extend the notion of singularity to cofinite Kleinian groups.  

In  \S\ref{sectionSelberg} we give an explicit form of the Selberg trace formula for cofinite Kleinian groups with finite-dimensional unitary representations (see Theorem \ref{T:Selberg}).
The new feature in the trace formula is a term of the form,
$$ \sum _{\alpha=1}^{\kappa}\frac{g(0)}{|\Gamma _{\alpha}:\Gamma '_{\alpha}|}\sum _{k=l_{\alpha}+1}^{\dim_\CC V}L(\Lambda_{\alpha},\psi_{k \alpha}). $$
The above term comes about from  \emph{regularity} at a cusp, and its value is computed using Kronecker's second limit formula ( see  \S\ref{sec.par}).
   
As an application of the spectral decomposition theorem we derive an identity  involving conjugacy relations of cuspidal elliptic elements.  This identity is used in the proof of the Selberg trace formula and to show that under certain conditions, the Selberg zeta function admits a meromorphic continuation (see Lemma \ref{lemCuspElip} for the identity).

For $\R(s)>1 $ the Selberg zeta-function $Z(s,\Gamma,\chi)$ is defined by the following product 
$$
Z(s,\Gamma,\chi) \df \prod_{ \{T_0 \} \in \mathcal{R}} ~ \prod_{j=1}^{ \dim_\CC V} \prod_{  \substack{ l,k \geq 0 \\  c(T,j,l,k)=1   } } \left( 1-\mathfrak{t}_{j} a(T_0)^{-2k} \overline{ a(T_0) ^{-2l}} N(T_0)^{-s - 1}    \right).  
$$
In \S\ref{secSZF} we introduce the various definitions and notations that are needed in order to define the zeta function, and   meromorphically\footnote{We give the meromorphic continuation for certain cases, and for others show that the zeta function is a rational root of a meromorphic function.}  continue    $Z(s,\Gamma,\chi)$ to $\R(s) \leq 1$ while computing its divisor. The main new difficulty is handling the contribution of the cuspidal elliptic elements to the topological (or trivial) zeros and poles of $Z(s,\Gamma,\chi).$  We show the following in \S\ref{secSZF}:
\begin{cor*}
Let $\Gamma = \PSL(2,\ZZ[-\frac{1}{2}+\frac{\sqrt{-3}}{2}]),$ and $\chi \equiv 1$ (the trivial representation). Then $Z(s,\Gamma,\chi)$ is not a meromorphic function\footnote{This is the first example that the author is aware of where the Selberg zeta-function is not meromorphic.} (it is the 6-th root of a meromorphic function).
\end{cor*}
In addition, the methods of \S\ref{secSZF} imply that the Selberg zeta-function of the Picard group is meromorphic:
\begin{cor*}
Let $\Gamma = \PSL(2,\ZZ[\sqrt{-1}]),$ and let $\chi \in \rep.$ Then $Z(s,\Gamma,\chi)$ is a meromorphic function.
\end{cor*}

All of the results of this paper appear in my PhD thesis \cite{Friedman} 

(\url{http://www.math.sunysb.edu/users/joshua/phdthesis}) 

with more details. 

I would like to thank my thesis advisor Professor Leon Takhtajan for the years he has spent guiding and teaching me.  I would also like to thank Professor J\"{u}rgen Elstrodt for reading this paper, and for many useful suggestions.  Special thanks  are due to Jay Jorgensen, Irwin Kra, Lee-Peng Teo, Alexei Venkov, and Peter Zograf for useful comments and suggestions.

\section{The Spectral Decomposition Theorem} \label{secSpectral}
Throughout this section $\Gamma $ is a cofinite Kleinian group and $\chi \in \rep.$

We start with the definition of a cusp. For every $\zeta \in \PP^1$  let  $\Gamma_\zeta $ denote the stabilizer subgroup of  $\zeta$ in $\Gamma, $ 
$$
\Gamma_{\zeta} \df \{~ \gamma \in \Gamma ~ | ~ \gamma \zeta = \zeta ~\},
$$ 
and  let $ \Gamma_\zeta^\prime $ be the maximal torsion-free parabolic subgroup of $\Gamma_\zeta $ (the maximal subgroup of $\Gamma_\zeta $ that does not contain elliptic elements). A point $\zeta \in \PP^1$ is called a \emph{cusp} of $\Gamma $ if  $\Gamma_\zeta^\prime$ is a free abelian group of rank two.  Two cusps $\zeta_1, \zeta_2$ are $\Gamma-$equivalent  if $\zeta_1 \in \Gamma \zeta_2,$ that is their $\Gamma-$orbits coincide.    

Every cofinite Kleinian group has finitely many equivalence classes of cusps, so we fix a set 
$ \{ \zeta_\alpha \}_{\alpha = 1}^{\kappa} $ of representatives of these equivalence classes.  For notational convenience we set $\Gamma_\alpha \df \Gamma_{\zeta_\alpha}, $ $\Gamma_\alpha^\prime \df \Gamma_{\zeta_\alpha}^\prime. $  
\begin{rem}
The possible values for the index of $[\Gamma_\alpha:\Gamma_\alpha^\prime]$ are 1,2,3,4, and 6.  See \cite[Theorem 2.1.8 (3), page 37]{Elstrodt}.
\end{rem}

For each cusp $\zeta_\alpha$ fix an element $B_\alpha \in \pc, $  a lattice $\Lambda_\alpha = \ZZ \oplus \ZZ \tau_\alpha,~\I(\tau_\alpha) > 0, $ and a root of unity $\epsilon_\alpha $ of order 1,2,3,4,or 6 with the following conditions being satisfied:

(1) $\zeta_{\alpha} = B_{\alpha}^{-1} \infty,$

(2) $$  B_{\alpha}\Gamma_{\alpha}^\prime B_{\alpha}^{-1} =  
\left\{ \,\left.\left(\begin{array}{cc}
1 &  b\\
0 &   1\end{array}\right)\,\right|\, b \in \Lambda_\alpha \,\right\},  $$

(3) $$  B_{\alpha}\Gamma_{\alpha} B_{\alpha}^{-1}= \left\{ \,\left.\left(\begin{array}{cc}
\epsilon & \epsilon b\\
0 & \epsilon^{-1}\end{array}\right)\,\right|\, b \in \Lambda_{\alpha},~\epsilon~\text{is some power of $\epsilon_\alpha$} \right\} /\{\pm I\}$$
The group $  B_{\alpha}\Gamma_{\alpha}^\prime B_{\alpha}^{-1} $ acts on  $\CC $ via the lattice $\Lambda_\alpha.$

See \cite[Theorem 2.1.8]{Elstrodt} for more details.

\begin{rem}
The group $ B_{\alpha}\Gamma_{\alpha}B_{\alpha}^{-1} $ also has a discrete action on $\CC.$  In fact it acts on every parallel copy of $\CC $ in $\HH$ by translations and rotations.  If $\gamma \in B_{\alpha}\Gamma_{\alpha}B_{\alpha}^{-1}$ and $P = z+rj \in \HH$ then $r=r(P) = r(\gamma P).$
\end{rem}

The notion of \emph{singularity} at a cusp\footnote{The group $\galp $ is insufficient for defining singularity of a cusp. We must use the full stabilizer subgroup $\gal.$} will play a prominent role in the proof of both the spectral decomposition theorem and the Selberg trace formula. 
For each cusp (representative)  $\zeta_\alpha$ of $\Gamma $ set the \emph{singular} space 
\beq \label{eqSingular}V_{\alpha} \df \{ v \in V \, | \, \chi(\gamma)v=v, \, \, \, \forall \gamma \in \Gamma_{\alpha}\,\}, \eeq  and the \emph{almost singular} space 
\beq \label{eqAlmostSingular} V_{\alpha}^\prime \df \{ v \in V \, | \, \chi(\gamma)v=v, \, \, \, \forall \gamma \in \Gamma_{\alpha}^\prime \,\},  \eeq  where $1 \leq \alpha \leq \kappa. $
\bd
A representation   $\chi \in \rep $  is  singular\footnote{See \cite{Roelcke} and \cite{Fischer} for the two-dimensional definition of singularity}  at the cusp $ \zeta_\alpha $ of $\Gamma $ iff the subspace 
$V_{\alpha}  \neq \{0 \}.$     
\ed
If a cusp is not \emph{singular} it is called \emph{regular}.
A representation $\chi $ is called \emph{singular} if it is singular  at least at one cusp, and \emph{regular}  otherwise. For each cusp $\zeta_\alpha$  set $k_\alpha = \dim_\CC V_\alpha, $   and
$$k(\Gamma,\chi) \df \sum_{\alpha = 1}^\kappa k_\alpha.  $$ 

We will need part (4) of the next lemma.  Its purpose is to help \emph{diagonalize} the not necessarily abelian group $\chi(\Gamma_\alpha).$
\begin{lem}
\label{L:mylemma} 
There exist $E_\alpha,R_\alpha,S_\alpha \in \Gamma_{\alpha} $ with the following properties:

(1)  $\Gamma_{\alpha}=\{\, E_\alpha^{k}R_\alpha^{i}S_\alpha^{j}\,|\,0\leq k<m_\alpha,\, i,j\in\ZZ\,\}.$
Here $R_\alpha,S_\alpha $ are parabolic elements with $B_\alpha R_\alpha B_\alpha^{-1} (P) = P+1$ and $B_\alpha S_\alpha B_\alpha^{-1} (P) = P+\tau_\alpha$ (here $\Lambda_\alpha = \ZZ \oplus \ZZ \tau_\alpha $)  for all $P \in \HH,$  and  $E_\alpha$ is elliptic of order $m_\alpha.$ 

(2)  $\Gamma_{\alpha}^\prime =\{\, R_\alpha^{i}S_\alpha^{j}\,|\, i,j\in\ZZ\,\}.$

(3) The elements $R_\alpha$ and $S_\alpha$ commute but the group $\Gamma_{\alpha}$ is not abelian when $m_\alpha>1$. 

(4)If in addition, $m_\alpha > 1, $  
then  $\chi(E_\alpha)$ maps $V_\alpha^\prime $ onto itself. Furthermore, there exists
a basis of $V_\alpha^\prime$ so that $\chi(E_\alpha)|_{V_\alpha^\prime}$ is diagonal. 
\end{lem}
\begin{proof}
(1), (2), and (3) readily follow from \cite[Theorem 2.1.8]{Elstrodt}. We prove (4):
set $E = E_\alpha, R=R_\alpha, S=S_\alpha.$ Since $ERE^{-1}$ and $R^{-1}$ are both parabolic and in  $\Gamma_{\alpha}^\prime $ it follows that 
 the $A \df E^{-1}R^{-1}E R \in\Gamma_{\alpha}^\prime $ is parabolic.    Since $A$ is parabolic, $A=R^{i}S^{j}$
for some $i,j\in\ZZ$ (applying (2)). By definition, the restriction $\chi(A)|_{V_\alpha^\prime }=I_{V_\alpha^\prime }.$ 
We have  
\[ \chi(E)\chi(R)=\chi(R)\chi(E)\chi(A). \] Next applying an arbitrary $v\in V_\alpha^\prime $ we obtain 
\beq
\chi(E)v = \chi(E)\chi(R)v = \chi(R)\chi(E)\chi(A)v = \chi(R)\chi(E)v 
\eeq
Thus $\chi(R)$ fixes $\chi(E)v$ and similarly $\chi(S)$ fixes
$\chi(E)v,$ so by definition  $\chi(E)v\in V_\alpha^\prime.$  Since $\chi(E)$ is
 unitary, its restriction to $V_\alpha^\prime$ is also unitary, hence
$V_\alpha^\prime$ has a diagonalizing basis. 
\end{proof}

Just as in the Fuchsian case, regularity of the representation implies discreteness of the spectrum of $\lp.$

\begin{lem} \label{lemRegResolvent}
Suppose $\Gamma $ is cofinite, and $\chi \in \rep. $    Then  for $\lambda \notin [0,\infty), $ the resolvent operator $R_\lambda = (\lp -\lambda)^{-1}$ is of Hilbert-Schmidt type iff $\chi $ is regular.
\end{lem}  
The proof follows from the Fourier series expansion of the resolvent kernel and part (4) of lemma \ref{L:mylemma}.  See \cite{Friedman} for more details.

In the \td case, the stabilizer subgroup of a cusp is a free abelian group of rank one contains only parabolic elements.    In our case $\gal$ contains a maximal rank-two  finite index subgroup $\galp$ and may also contain  \emph{cuspidal elliptic} elements.  Surprisingly, by Lemma~\ref{lemRegResolvent}, the cuspidal elliptic elements can affect the spectral\footnote{Suppose $\Gamma$ has only one cusp at $\infty.$ If $\chi$ is a character that is trivial on $\gip,$ then $\lp$ has a continuous spectrum iff $\chi$ is trivial on all cuspidal elliptic elements of $\gi.$  } properties of $\lp.$  This is seen from Lemma~\ref{lemRegResolvent}. 

For $\cusps $ fix a basis $ \{ v_l(\alpha) \}_{l=1}^{k_\alpha} $ for $V_\alpha. $ Suppose that   $P \in \HH,\, \R(s)>1, \cusps, $ and $\sing,$ we define the  \emph{Eisenstein series}  by 
$$
E_{\alpha l}(P,s)   \df  E(P,s,\alpha ,l,\Gamma,\chi) \df  \sum_{M \in \Gamma_\alpha \setminus \Gamma } \left(r(B_{\alpha}MP)
\right)^{1+s}\chi(M)^{*}v_l(\alpha). $$

The series $E_{\alpha l}(P,s) $ converges uniformly and absolutely on compact subsets of  $\{ \R(s)> 1 \} \times \HH $ to a $\chi-$automorphic function  that  satisfies
$$ \lp E(~\cdot~,s,\alpha ,v) = \lambda E(~\cdot~,s,\alpha ,v),  $$  and admits a meromorphic continuation to the whole complex plane.  

The meromorphic continuation of $E_{\alpha l}(P,s) $ is necessary for the proof of the  spectral decomposition theorem, and is highly non-trivial.  Fortunately, there are several well-known methods available,  \cite{Faddeev}, \cite{Selberg2}, \cite{Selberg3}, \cite{Verdi}, and \cite{Langlands}.  In \cite{Elstrodt} an adaptation of the methods in \cite{Verdi} is used to prove the \thd scalar case, and a similar adaptation works\footnote{Alternatively, Faddeev's method can also be used by to express the meromorphic continuation of $E_{\alpha l}(P,s) $ by adapting \cite[Chapter 2 and 3]{Venkov}.} for the vector case (see \cite{Friedman}). Adapting \cite{Elstrodt}[Chapter 6] gives us:
\begin{thm}\label{T:spectral}(The Spectral Decomposition Theorem)
Let $f \in \hs.$   Then $f$ has an expansion of the form
\begin{equation}\label{E:spectral}
 f(P) = \sum_{m \in \D} \left< f,e_m \right> e_m(P) + 
\frac{1}{4\pi}\sum_{\alpha = 1}^\kappa \sum_{l = 1}^{k_\alpha}
\frac{\left[ \Gamma_\alpha:\Gamma_\alpha^\prime \right]}{\left| \Lambda_\alpha \right|}
\int_{\RR} \left<f,E_{\alpha l}(\,\cdot \,,it) \right> E_{\alpha l}(P,it) \, dt.
\end{equation}
 The sum and integrals converge in the Hilbert space $\hs.$ Here $\D$ is an indexing set of the eigenfunctions $e_m$ of $\lp$ with corresponding eigenvalues $\lambda_m,$   $E_{\alpha l}(P,s)$ are the Eisenstein series associated to the  singular cusps of $\Gamma$,   $k_\alpha = \dim_\CC V_\alpha $, and $\left| \Lambda_\alpha \right| $ is the Euclidean area of a fundamental domain for the lattice  $ \Lambda_\alpha. $ If a cusp is regular it is omitted from the sum in \eqref{E:spectral}. 
\end{thm} 

\section{The Selberg Trace Formula} \label{sectionSelberg}
\begin{thm}{(Selberg Trace Formula)} \label{T:Selberg}
Let $\Gamma $ be a cofinite Kleinian group,  $\chi \in \rep,~h $ be  a holomorphic function on  $ \{ s \in \CC \, | \, |\I(s)| < 2+ \delta \}$ for some $\delta > 0,$ satisfying $ h(1+z^2) = O( 1+|z|^2)^{3/2 - \epsilon}) $ as  $|z| \ra \infty,$ and let
$$ g(x) = \frac{1}{2\pi} \int_{\RR} h(1+t^2)e^{-itx}\,dt. $$  Then
\begin{multline}
\sum _{m \in \D }h(\lambda _{m}) 
= \frac{\vol \left( \Gamma \setminus \HH \right)}{4\pi ^{2}}\dim_\CC V \int _{\RR }h(1+t^{2})t^{2}\, dt  \\ + 
 \sum_{ \{R \} \text{\emph{nce}}}\frac{\tr _{V}\chi (R) g(0)\log N(T_{0})}{4|\ren |\sin ^{2}(\frac{\pi k}{m(R)})}  + 
\sum_{\{ T \} \text{\emph{lox}} } \frac{\tr _{V}\chi (T) g(\log N(T))}{\oen |a(T)-a(T)^{-1}|^{2}}\log N(T_{0}) \\ - 
\frac{\tr (\smat (0))h(1)}{4}+\frac{1}{4\pi }\int _{\RR }h(1+t^{2})\frac{\phi'}{\phi }(it)\, dt  \\ +
\sum_{\alpha=1}^{\kappa}\sum _{k=1}^{e_{\alpha}}\frac{\tr _{V}\chi (g_{\alpha k})}{|C(g_{\alpha k})|}\left(g(0)c_{\alpha k}+d_{\alpha k}\int _{0}^{\infty }g(x)\frac{\sinh x}{\cosh x-1+\alpha _{\alpha k}}\, dx\right)  \\ + 
\sum _{\alpha =1}^{\kappa}\left(\frac{l_{\alpha}}{|\Gamma _{\alpha}:\Gamma '_{\alpha}|}\left(\frac{h(1)}{4}+g(0) (\frac{1}{2}k_{\Lambda_{\alpha}}-\gamma )- \frac{1}{2\pi }\int _{\RR }h(1+t^{2})\frac{\Gamma '}{\Gamma }(1+it)\, dt\right)\right)  \\ +
\sum _{\alpha=1}^{\kappa}\frac{g(0)}{|\Gamma _{\alpha}:\Gamma '_{\alpha}|}\sum _{k=l_{\alpha}+1}^{\dim_\CC V}L(\Lambda_{\alpha},\psi_{k \alpha}).
\end{multline}
\end{thm}
Here, $ \{ \lambda _{m} \}_ {m \in \D } $ are the eigenvalues of $\lp$ counted with multiplicity.  Following \cite{Elstrodt} section 5.2,  the summation with respect to $\{R \}_\text{nce} $ extends over the finitely many $\Gamma-$conjugacy classes of the non cuspidal elliptic elements (elliptic elements that do not fix a cusp) $R \in \Gamma,$ and for such a class $N(T_0)$ is the minimal norm of a hyperbolic or loxodromic element of the centralizer $\mathcal{C}(R).$  The element $R$ is understood to be a $k-$th power of a primitive non cuspidal elliptic element $R_0 \in \mathcal{C}(R)$ describing a hyperbolic rotation around the fixed axis of $R$ with minimal rotation angle $\frac{2 \pi}{m(R)}.$  Further, $\ren$ is the maximal finite subgroup contained in $  \mathcal{C}(R).$  The summation with respect to $\{ T \}_\text{lox}$ extends over the  $\Gamma-$conjugacy classes of hyperbolic or loxodromic elements of $\Gamma,$ $T_0$ denotes a primitive hyperbolic or loxodromic element associated with $T.$  The element $T$ is conjugate in $\pc$ to the transformation  described by the diagonal matrix with diagonal entries $a(T), a(T)^{-1}$ with $|a(T)| > 1, $ and $N(T) =   |a(T)|^2. $ For $s \in \CC,$ $\smat(s)$ is  a $k(\Gamma,\chi) \times  k(\Gamma,\chi)$ matrix-valued meromorphic function, called the \emph{scattering matrix}\footnote{The definition of the scattering matrix uses the \emph{Fourier expansion of a $\chi$-automorphic at a singular cusp} which though not difficult, is long and combersome.  The corresponding \tdv case defintion is almost identical to our case and we refer the reader to  \cite[Chapters 2 and 3]{Venkov} and \cite{Friedman}. }  of $\lp,$ and $\phi(s) = \det \smat(s).$  Furthermore $c_{\alpha k}, g_{\alpha k},$ and $d_{\alpha k}$ are constants depending on $\Gamma$ which will be determined in the case of $\Gamma $ having only one cusp at $\infty. $  The remaining notation will be defined in  \S \ref{sec.par}\footnote{Please note that there is a typographical error in the loxodromic and non cuspidal elliptic terms in \cite{Elstrodt} Theorem 6.5.1; both terms are missing a factor of $\frac{1}{4 \pi}. $}.

\section*{The Spectral Trace}

In order to make this paper self-contained, we outline the (standard) proof of the spectral side of the Selberg trace formula.  

For $P = z+rj,~P'=z'+r'j \in \HH$ set 
$$ \delta(P,P') \df  \frac{|z-z'|^{2}+r^{2}+ r'^{2}}{2rr'}. $$ It follows that $\delta(P,P') = \cosh(d(P,P'))$, where $d$ denotes the hyperbolic distance  in $\HH.$  Next, for $k \in \scz([1,\infty))$ a Schwartz-class function, define 
$$K(P,Q) = k(\delta(P,Q)),~\text{and}~ K_\Gamma(P,Q) \df \sum_{\gamma\in\Gamma} \chi(\gamma)K(P,\gamma Q). $$ 

The series above converges absolutely and uniformly on compact subsets of
$\HH \times \HH$, and is the kernel of a bounded operator $\mathcal{K} : \hs \mapsto \hs.$  The Selberg trace formula is essentially\footnote{We say ``essentially'' because $\K$ is not of trace class.  Selberg's  procedure is used to define and compute the \emph{regularized trace.}} the trace of $\mathcal{K}$ evaluated in two different ways: the  first using spectral theory, and the second as an explicit integral.

The function $h$ that appears in the Selberg trace formula is the  Selberg--Harish-Chandra transform\footnote{If $f:\HH \mapsto V$ is a smooth function satisfying $\lp f = \lambda f$, then $\mathcal{K}f= h(\lambda)f.$  That is $f$ is an eigenfunction of $\mathcal{K}$ with an eigenvalue that depends only on $\lambda.$ } of $k,$ defined as follows:  
\beq \label{eqSHC}
h(\lambda) = h(1-s^2) \df \frac{\pi}{s}
\int_{1}^{\infty}k\left(\frac{1}{2}\left(t+\frac{1}{t}\right)\right)
(t^{s}-t^{-s})\left(t-\frac{1}{t}\right)\,\frac{dt}{t},~~\lambda = 1-s^2.
\eeq

The first step in evaluating the spectral trace is to compute the spectral expansion of  $\K.$  For $v,w \in V $ let $ v \otimes \overline{w}$ be a linear operator in $V$ defined by    
$ v \otimes \overline{w}(x) = <x,w>v.  $
An immediate application of the spectral decomposition theorem (Theorem \ref{T:spectral}) and the Selberg--Harish-Chandra (Equation \ref{eqSHC}) transform  gives us (see \cite[Equation 6.4.10, page 278]{Elstrodt}),
\begin{lem}
Let  $k \in \scz $ and  $h:\CC\ra\CC$ be the Selberg--Harish-Chandra  Transform
of $k.$ Then with $K_\Gamma$ defined above we have 
\begin{multline}
\label{E:kernel expansion}
K_\Gamma(P,Q)   = \sum_{m \in \D} h(\lambda_m)e_m(P) \otimes 
\overline{e_m(Q)}   \\ +
\frac{1}{4\pi}\sum_{\alpha = 1}^\kappa \sum_{l = 1}^{k_\alpha}
\frac{\left[ \Gamma_\alpha:\Gamma_\alpha^\prime \right]}
{\left| \Lambda_\alpha \right|}
\int_{\RR} h \left( 1+t^2 \right) E_{\alpha l}(P,it) \otimes \overline{E_{\alpha l}(Q,it)} \, dt.
\end{multline}
The sum and integrals converge absolutely and uniformly on compact subsets 
of $\HH \times \HH$.
\end{lem}
Next, we split up $K_\Gamma $ as a sum of two kernels. The first kernel 
$$
H_\Gamma(P,Q) =  \frac{1}{4\pi}\sum_{\alpha = 1}^\kappa  
\sum_{l = 1}^{k_\alpha}
\frac{\left[ \Gamma_\alpha:\Gamma_\alpha^\prime \right]}
{\left| \Lambda_\alpha \right|}
\int_{\RR} h \left( 1+t^2 \right) E_{\alpha l}(P,it) 
\otimes \overline{E_{\alpha l}(Q,it)} \, dt,
$$ is not of Hilbert-Schmidt class, while the second kernel $$
L_\Gamma(P,Q) =  \sum_{m \in \D} h(\lambda_m)e_m(P) \otimes 
\overline{e_m(Q)}
$$  
is of trace class since $h$ decays sufficiently fast\footnote{This follows from the fact that $k$ is a Schwartz class function.}.

In order to define the regularized trace of $\K$ we need an explicit description of a fundamental domain for $\Gamma,$ particularly  in the cusp sectors. Suppose $Y>0$ is sufficiently large. Then for all $A > Y$ there exist a compact set $\F_{A}\subset\HH$
such that \[
\F \df \F_{A}\cup\F_{1}(A)\cup\cdots\cup\F_{\kappa}(A)\]
 is a fundamental domain for $\Gamma.$  The sets $\F_{\alpha}(A) $ are cusp sectors of each cusp of $\Gamma $ (see \cite{Elstrodt} Proposition 2.3.9).  We can now compute the  truncated trace.
\bl \label{lemSingCancel}
\begin{multline} 
\lim_{A \ra \infty} \left( \int_{\F_A} \tr_V (K_{\Gamma}(P,P))  \,dv(P)
-  \int_{\F_A} \tr_V (H_{\Gamma}(P,P))  \,dv(P)  \right) \\
=  \int_\F \tr_V (L_{\Gamma}(P,P))\,dv(P) = 
\sum_{m \in \D} h(\lambda_m)  < \infty.
\end{multline} 
The infinite sum is absolutely convergent.  In particular any divergent terms (as $A \ra \infty$) are canceled out.
\el
\begin{proof}
Its not hard to see that the restriction of the resolvent kernel (of $\lp$) to the span of the set of eigenfunctions of $\lp$ is a Hilbert-Schmidt kernel (see \cite[Theorem 4.5.2]{Elstrodt}).  Thus it follows that 
$$\sum_{m \in \D} \frac{1}{1+|\lambda|^2}  < \infty.  $$  Since $k$ is of rapid decay, so is $h$ (by Equation \ref{eqSHC}).  Thus 
$$\sum_{m \in \D} h(\lambda_m)  < \infty.  $$   
\end{proof}
We will apply the above lemma in \S\ref{sec.par} to compute some  identities involving cuspidal elliptic elements.

Next, applying the vector form of  the Maa\ss-Selberg  relations (see \cite[Pages 67-70]{Venkov} and \cite[Page 305]{Elstrodt} for the similar \tdv and \thds cases)  we can   compute the spectral truncated trace:
\begin{lem} \label{L:SpecTrace}
\begin{multline}
\int_{\F_A} \tr_V (H_{\Gamma}(P,P))  \,dv(P) = \\ 
g(0)k(\Gamma,\chi) \log(A) - \frac{1}{4\pi} \int_{\RR}\lds(it)
h(1+t^2)\,dt + \frac{h(1)\tr \smat(0)}{4} +  
\lto_{A \ra \infty}.
\end{multline}
The integral on the right hand side converges absolutely.
\end{lem}

Here $\smat(s) $ is the scattering matrix, and $\phi(s) = \det \smat(s).$  

\section{The Explicit Trace}
In this section our main goal is to give an explicit formula for $\int_{\F_A} \tr_V (K_{\Gamma}(P,P))  \,dv(P).  $  By Lemma \eqref{L:SpecTrace} and Lemma \eqref{lemSingCancel}, there is a constant $C_k$ so that  
\beq 
\int_{\F_A} \tr_V (K_{\Gamma}(P,P))  \,dv(P) = g(0)k(\Gamma,\chi) \log(A) + C_k +  \lto_{A \ra \infty}.
\eeq
We determine $C_k$ by explicitly integrating $\int_{\F_A} \tr_V (K_{\Gamma}(P,P))  \,dv(P), $ following Selberg's original method.

We  decompose 
$$K_\Gamma(P,Q) = \sum_{\gamma \in \Gamma} \chi(\gamma)
K(P,\gamma Q)$$  into various sub-sums.  Depending on their \emph{type.} The types are as follows,   ``id'' is the identity,  ``par''  are the parabolic elements, ``ce'' are the cuspidal elliptic elements\footnote{These elliptical elements share a common fixed point in $\PP$ with some parabolic element in $\Gamma.$}, ``nce'' are the non-cuspidal elliptical elements, ``lox'' are the hyperbolic and loxodromic elements, and ``cusp'' $=$ ``par'' $\union$ ``ce''.
For each $ S \in \{\ID,\PAR,\CE,\NCE,\LOX,\CUSP\}$ set 
$$
K_{\Gamma}^{S}(P,Q) = \sum_{\gamma \in \Gamma^S} \chi(\gamma)
K(P,\gamma Q).
$$
Here $\Gamma^S$ denotes the subset of $\Gamma $ consisting of elements of type $S.$  

Following \cite{Elstrodt} section 5.2 and Theorem 6.5.1 we have,
\begin{lem} \label{lemNonCuspidal}
\begin{multline}
\int_{\F_A} \tr_V \left( K_{\Gamma}^{\ID}+ K_{\Gamma}^{\NCE}
+  K_{\Gamma}^{\LOX}  \right)(P,P)  \,dv(P)  \\
= \frac{ \vol (\Gamma \setminus \HH) \dim_\CC (V) } { 4 \pi^2 } 
\int_{\RR} h(1+t^2)t^2 \, dt 
+ \sum_{ \{ R \}_{\NCE} } \frac{ \tr(\chi(R)) g(0) \log N(T_0) }{ 4 |\E(R)| \sin^2 \left(\frac{ \pi k }{ m(R) } \right) } 
\\ + \sum_{ \{ T \}_{\LOX}} \frac{\tr(\chi(T)) g( \log N(T))
\log N(T_0)}{|\E(T)| |a(T) - a(T)^{-1}|^2}  + o(1) ~\text{as $A \ra \infty.$}
\end{multline}
\end{lem}
The notations above were defined in \S\ref{sectionSelberg}.

\section*{The Cuspidal Elliptic Elements}
Our next immediate goal is to evaluate  
$$ \int_{\F_A} \tr_V  K_{\Gamma}^{\CE}(P,P)  \,dv(P).  $$

For notational simplicity, we will adopt the following assumption   from this point on until the end of this paper.
\begin{as} \label{asOne}
The Kleinian group $\Gamma$ has only one class of cusps at $\zeta = \infty \in \PP,$ and $\chi \in \rep.$
\end{as}

Denote by $ \cuspi $  set of elements of $\Gamma$ 
which are $\Gamma$-conjugate to an element of $ \Gamma_\infty
\setminus \Gamma_\infty^\prime = \{ \gamma \in \gi~|~\gamma~\text{is not parabolic nor the identity element}~\}. $  
We fix representatives  of conjugacy classes of  $\cuspi,~
g_{1}, \dots , g_{d}\footnote{There are only finitely many distinct conjugacy classes of elliptic elements in a cofinite Kleinian group.} $ that have the form 
\beq g_{i} =  \left(\begin{array}{cc}
 \epsilon_{i} & \epsilon_{i} \omega_{i} \\
 0 & \left( \epsilon_{i}  \right)^{-1}\end{array}
\right).
\eeq
Let $\mC (g) $ denote the centralizer in $\Gamma $ of an element $g \in \cuspi.$  In addition, let   $\{ p_i,\infty \} $ be the set of fixed points in $\cinf$ of the element $g_i.$  Since $g_i $ is a cuspidal elliptic element it follows that  $p_i $ is a cusp of $\Gamma $ (see \cite{Elstrodt} page 52).  Hence by Assumption~\ref{asOne} there is an element $\gamma_i \in \Gamma $ with $\gamma_i \infty = p_i. $  Suppose that $c_i $ is the lower left hand (matrix) entry of $\gamma_i.$ Then we have (see \cite[Pages 302-304]{Elstrodt}),
\bl \label{lemCuspidalElliptic}
\begin{multline}
\int_{\F_A} \tr K_\Gamma^{\CE}(P,P)~ dv(P) \\ = \sum_{i=1}^d \frac{\tr \chi(g_i)}{|\mC(g_i)|}\left(  \frac{2 g(0) (\log|c_i| + \log A)}{|1-\epsilon_i^2|^2}+\frac{1}{|1-\epsilon_i^2|^2} \int_0^\infty g(x) \frac{\sinh x}{\cosh x - 1 +\frac{|1-\epsilon_i^2|^2}{2} }~dx \right) \\ + o(1) ~\text{as $A \ra \infty.$}
\end{multline}
\el

\section{The Contribution of the Parabolic Elements to the Selberg Trace Formula} \label{sec.par}
This section contains the new features of the Selberg trace formula that are not present in the \tdv and \thds cases.  We remind the reader that Assumption~\ref{asOne} is in effect.

Our main goal for this section is to evaluate 
\beq \int_{\F_A} \tr_V  K_{\Gamma}^{\PAR}(P,P)  \,dv(P). \label{eqParInt} \eeq
The computation of \eqref{eqParInt} can be split up into two parts.  The first part is the evaluation of 
$$
\int_{\F_A} \tr_V  K_{\Gamma}^{\PAR}(P,P)\mathbf{v}  \,dv(P)
$$
for $ \mathbf{v} \in V_\infty^\prime,$ and is an immediate extension of the \thds case.  The other part is the computation of $$
\int_{\F_A} \tr_V  K_{\Gamma}^{\PAR}(P,P) \mathbf{w}  \,dv(P)
$$
for $\mathbf{w} \in \left( V_\infty^\prime  \right)^\perp, $ and requires us to study a particular lattice sum.

\section*{Lattice Sums}
Let $\Lambda =  \ZZ \oplus \ZZ \tau \subset \CC$ be a lattice with $\I(\tau) > 0. $  A (lattice) character  $ \psi $ of $\Lambda $ is a one-dimensional unitary representation of $\Lambda. $
\bd
For  $x > 0 $ set 
$$ Z(x,\Lambda , \psi) \df \sum_{ \substack{ \mu \in \Lambda \\  | \mu |^2 \leq x \\ \mu \neq 0 } } \frac{ \psi(\mu)}{ |\mu |^2},    $$   
and when the limit exists 
$$  L(\Lambda , \psi ) \df  \lim_{x\ra \infty} Z(x,\Lambda , \psi). $$
\ed

\bp \label{P:latticesum} Let $ \psi $ be a  character of $\Lambda $. 

(1) If $\psi = \text{id}$, the trivial character, then 
$$ Z(x,\Lambda , \psi) = \frac{\pi}{|\Lambda |}  (\log x + \kappa_\Lambda ) + O \left( x^{ - \frac{1}{2}} \right) \space  \text{ \space \space \space as $x \ra \infty $.}  $$

(2) If $\psi \neq \text{id}$  then $ \lim_{x\ra \infty} Z(x,\Lambda , \psi) $ exists and 
$$ Z(x,\Lambda , \psi) = L(\Lambda , \psi ) + O \left( x^{- \frac{1}{2}} \right) \space  \text{ \space \space \space as $x \ra \infty $.}  $$
\ep
Here $ \kappa_\Lambda $ is an analogue of the Euler constant for the lattice $\Lambda.$ 
\pf
(1) is proved in \cite[Lemma 5.2, page 298]{Elstrodt}. (2) follows from a long, but straight forward summation by parts argument.  The details are worked out in \cite{Friedman}.
\epf
Our next goal is to determine the value of $ L(\Lambda , \psi ).$  We will need  Kronecker's second limit formula.
Let $u,v$ be real numbers which are not both  integers, and let $ \tau = x + i y,~y>0.$   For $ \R(s) > 1 $ set 
$$  E_{u,v}(\tau, s) \df   \psum_{m,n=-\infty}^\infty  e^{ 2\pi i (mu+nv)} \frac{y^s}{|m \tau + n|^{2s}}. 
$$
The sum is defined over all integers $n,m,$  and the prime in the sum means to leave out  $(0,0).$
The series converges  uniformly and absolutely on compact subsets of $ \R(s) > 1. $   We have  (\cite{Siegel} or  \cite[page 276]{Lang}\footnote{There appears to be a typographical error in the definition of the Siegel function on page 276 of the second edition of \cite{Lang}.  The correct definition appears on page 262.})
\begin{lem} \label{P:kron}
The function $ E_{u,v} (\tau, s) $ can be continued to an entire function of $s \in \CC,$ and one has   \beq  E_{u,v} (\tau, 1) = - 2 \pi \log | g_{-v,u} (\tau) |, \eeq 
where $g_{a_1,a_2} $ is the Siegel function, \beq
g_{a_1,a_2}(\tau) = -q_{\tau}^{(1/2)\textbf{B}_{2} (a_1)} e^{ 2 \pi i a_2 (a_1-1)/2 }(1-q_z) \prod_{n=1}^{\infty} (1 - q_{\tau}^n q_z)(1 - q_{\tau}^n /q_z),
\eeq
$ \textbf{B}_{2} (X) = X^2 -X + 1/6, $ $ q_{\tau} = e^{ 2 \pi i \tau},$  $ q_z =  e^{ 2 \pi i z}, $ and $z = a_1\tau+a_2. $
\end{lem}
Finally we can evaluate $L( \Lambda, \psi ).$
\bp Let $\Lambda =  \ZZ \oplus \ZZ \tau \subset \CC$ be a lattice with $\I(\tau) > 0, ~ \psi $ a character of  $\Lambda,$  and $u,v \in \RR $ are not both  integers satisfying  
$  \psi(1) =  e^{2 \pi i u} ~~\text{and}  $
$ \psi(\tau) = e^{2 \pi i v}. $
Then 
\beq L( \Lambda, \psi ) = \frac{-2 \pi}{y} \log \left| g_{-v, u} \left( \tau \right) \right|.      \eeq
\ep
\pf
Using a summation by parts argument (see \cite{Friedman}) one can show that for $s \in [1,2]$ the function\footnote{We introduce new notation because for the value $s=1$ a priori we do not know that analytic continuation of $ E_{u,v} (\tau, s) $ agrees with $f(1).$ Once we prove $f$ is continuous the new notation will be redundant.}
$$  f(s) \df  \lim_{x \ra \infty} \psum_{|m \tau + n|^{2} < x}  e^{ 2\pi i (mu+nv)} \frac{y^s}{|m \tau + n|^{2s}}. 
$$ converges uniformly (though not absolutely) on the interval $[1,2].$  Thus  $f(s)$ is continuous on $[1,2]$ and agrees with $ E_{u,v}(\tau, s) $ on $(1,2],$ hence we can pass to the limit as $s \ra 1^+$ and conclude that 
$$L( \Lambda, \psi ) =\frac{1}{y}E_{u,v}(\tau, 1).  $$
\epf

\section*{Evaluation of Integral \ref{eqParInt} }

Let  $\pg$ be a fundamental domain for the action\footnote{See \S\ref{secSpectral} for more details on the action.}  of  $\Gamma_{\infty} $ on $\CC $,  
$$\widetilde{\pg} \df  \{ (z,r) \in \HH \, | \, z \in   \pg \, \},$$ and $$ \pg_{A} \df    \{ (z,r) \in \HH \, | \, z \in   \pg, \, r \leq A \,  \}. $$
It follows that 
$\widetilde{\pg} $ is a fundamental domain for the action of  $ \Gamma_{\infty} $  on $\HH $. 

Recall that $ \gip $ is canonically isomorphic to a lattice $  \Lambda_{\infty}. $  For $\mu \in \Lambda_{\infty}$ let $\widehat{\mu}$ denote the corresponding parabolic element in $\gip. $ We will need the following (see \cite[Pages 300-301]{Elstrodt})
\bl \label{lemParIntToSum}
\beq   \int_{\F_A} \tr K_{\Gamma}^{\PAR}(P,P) \, dv(P) =  \psum_{\mu \in \Lambda_{\infty}} \tr( \chi ( \widehat{\mu} ) ) \int_{ \pg_{A}} K(P, \widehat{\mu} P) \,dv(P) + \lto_{A \ra \infty}. \eeq
\el
  
Since $\gip $ is an abelian group, $\chi $ restricted to $\gip $ can be diagonalized.  In other words, there exist lattice characters $ \{ \psi_l \}_{l = 1\dots n} $ so that 
\beq \label{eqTraceChi}
\tr \chi |_{\gip} = \sum_{l=1}^n \psi_l.  
\eeq  Thus it suffices to consider lattice characters instead of unitary representations.
\bl \label{lemParMain} Let $\psi $ be a lattice character of $\Lambda_\infty. $  Then 

(1) For $ \psi = \text{id}, $   
\begin{multline*}
 \psum_{ \mu \in \Lambda_\infty } \psi (\mu)  \int_{\pg_{A} }K(P,  \widehat{\mu} P) =  \\
\frac{1}{[ \Gamma_{\infty}:\Gamma_{\infty}^{\prime} ]} 
\left(  g(0) \log A + \frac{h(1)}{4} + 
g(0) \left( \frac{ \eta_{\infty}}{2} 
- \gamma \right) - 
\frac{1}{2\pi} 
\int_{\RR} h(1+t^2) \frac{\Gamma'}{\Gamma}(1+it) \,dt   \right) \\ + o(1) ~\text{as $A \ra \infty $}.
\end{multline*}

(2) For  $ \psi \neq \text{id}, $ 
$$
 \psum_{ \mu \in \Lambda_\infty} \psi (\mu)  \int_{\pg_{A} }K(P,  \widehat{\mu} P) = 
\frac{g(0)}{[\Gamma_{\infty}:\Gamma_{\infty}^{\prime}]} L(\Lambda, \psi ) +o(1) ~\text{as $A \ra \infty $}.
$$
\el
\pf
(1) is proved in \cite[pages 300-302]{Elstrodt}.

The proof of (2) is  a modification of (1).   Let $C_\infty = \frac{ | \Lambda_\infty  |}{ [ \Gamma_{\infty}:\Gamma_{\infty}^{\prime} ]}, $ then  by the definition of the action of $\widehat{\mu}$ on $\HH, $
\begin{multline} \label{eq393}
 \psum_{ \mu \in \Lambda_\infty } \psi (\mu)  \int_{\pg_{A} }K(P,  \widehat{\mu} P) = 
C_\infty \psum_{ \mu \in \Lambda_\infty } \psi( \mu ) \int_0^A k \left( \frac{ | \mu |^2 }{2 r^2} + 1 \right) \, \frac{dr}{r^3}.
\\ = C_\infty \psum_{ \mu \in \Lambda_\infty } \frac{ \psi (\mu) }{ | \mu |^2 } \int_{ \frac{| \mu |^2 }{2A^2}}^\infty k(u+1) \,du. 
\end{multline}

Since $$ Z(x,\Lambda , \psi) = \sum_{ \substack{ \mu \in \Lambda_\infty  \\  | \mu |^2 \leq x \\ \mu \neq 0 } } \frac{ \psi(\mu)}{ |\mu |^2}, $$  using summation by parts, we can rewrite
\eqref{eq393} as 
\beq
C_\infty \int_0^\infty k(u+1) Z(2A^2u,\Lambda, \psi )\,du.
\eeq
Next we apply Proposition \ref{P:latticesum} to  obtain 
\beq
C_\infty \int_0^\infty k(u+1) Z(2A^2u,\Lambda, \psi )\,du = C_\infty \int_0^\infty k(u+1) \left( L(\Lambda, \psi) + O \left( (2A^{2}u)^{-1/2} \right) \right) \,du. 
\eeq
Now we  show that the resulting error term is $ o(1). $  Since $ k(u+1) = O \left( (1+u)^{-4}    \right)  $  ($k$ is a rapid decay function)
$$ \left|  \int_0^\infty k(u+1) O \left( (2A^{2}u)^{-1/2} \right) \,du \right| \leq 
\frac{D}{A} \int_0^\infty (1+u)^{-4} u^{-1/2} \,du = o(1) \,\, \text{as} \,\, A \ra \infty 
 $$ for some $D > 0.$
To complete the proof note that $ \int_0^\infty k(u+1) \,du = g(0). $
\epf
Finally we can evaluate \eqref{eqParInt}:
\bl \label{lemParabolic}
\begin{multline}
\int_{\F_A} \tr K_{\Gamma}^{\PAR}(P,P) \, dv(P) =  \\
\frac{l_\infty }
{| \Gamma_{\infty}: \Gamma_{\infty}^{\prime} |} 
\left(  g(0) \log A + \frac{h(1)}{4} + 
g(0) \left( \frac{ \eta_{\infty}}{2} 
- \gamma \right) - 
\frac{1}{2\pi} 
\int_{\RR} h(1+t^2) \frac{\Gamma'}{\Gamma}(1+it) \,dt   \right)  \\ + 
 \frac{g(0)}{| \Gamma_{\infty}: \Gamma_{\infty}^{\prime} |} \sum_{ l = l_\infty + 1 }^{n} L(\Lambda_{\infty}, \psi_{ l}  ). 
\end{multline}
Here $n=\dim_\CC V, $ $\psi_{ l}$ are the lattice characters associated to the lattice $\Lambda_\infty, $ $l_\infty = \dim_\CC V_\infty^\prime, $ and  $\eta_\infty $ is the analogue of the Euler constant for the lattice $\Lambda_\infty. $ 
\el
\pf
The proof follows immediately from Lemma~\ref{lemParIntToSum}, Equation~\ref{eqTraceChi}, and Lemma~\ref{lemParMain}.
\epf

We have evaluated the truncated trace of $\mathcal{K}$ explicitly as an integral, and by using spectral theory.  Notice that as $A \ra \infty $ the integral over the parabolic sum (Lemma~\ref{lemParabolic}) has a divergent term.  So does the corresponding cuspidal elliptic integral ( Lemma~\ref{lemCuspidalElliptic}). By Lemma \ref{lemSingCancel} the divergent terms \emph{must} equal the divergent term of the spectral (truncated) trace (Lemma~\ref{L:SpecTrace}).  It follows that
$$2 g(0) \log A \sum_{i=1}^d \frac{\tr \chi(g_i)}{|\mC(g_i)||1-\epsilon_i^2|^2} +  
 g(0) \log A \frac{l_\infty }
{| \Gamma_{\infty}: \Gamma_{\infty}^{\prime} |}  - g(0) k_\infty \log A = 0. $$  
By choosing a suitable $k$ so that $g(0) \neq 0 $ we obtain
\bl \label{lemCuspElip}
$$
2  \sum_{i=1}^d \frac{\tr \chi(g_i)}{|\mC(g_i)||1-\epsilon_i^2|^2} +  
  \frac{l_\infty }
{| \Gamma_{\infty}: \Gamma_{\infty}^{\prime} |}  =  k_\infty. 
$$
\el
The formula\footnote{A similar formula is valid for the general case of $\kappa$-many cusps. } above is an  application of spectral theory to the group relations of a cofinite hyperbolic three-orbifold\footnote{Notice that all of the terms above are defined simply in terms of group relations. }.   We will use the above lemma  to give a meromorphic continuation of the Selberg zeta-function.

\section*{Completion of the proof of the Selberg Trace Formula}
The Selberg trace formula now follows: combine    Lemma~\ref{lemParabolic}, Lemma~\ref{lemCuspidalElliptic}, Lemma~\ref{lemCuspElip}, Lemma~\ref{lemNonCuspidal}, Lemma~\ref{L:SpecTrace}, and Lemma~\ref{lemSingCancel}.  Note that the divergent terms all cancel by Lemma~\ref{lemSingCancel} (or we can use  Lemma~\ref{lemCuspElip}).  Finally take the limit as $A \ra \infty.$  See \cite[Section 6.5]{Elstrodt} for more details on combining the lemmas above.

\section{The Selberg Zeta Function} \label{secSZF}
In this section we define the Selberg zeta-function $Z(s,\Gamma,\chi)$    for cofinite Kleinian groups with finite-dimensional unitary representations, in the right half-plane $\Re(s) > 1.$  We then  evaluate the logarithmic derivative of  $Z(s,\Gamma,\chi)$ and show that $Z(s,\Gamma,\chi)$ admits a meromorphic continuation, subject to some technical assumptions concerning the stabilizer subgroup $\gi.$

\section*{The Definition of the Selberg Zeta-Function and its Motivation}

In the celebrated paper \cite{Selberg1} Selberg first defined what is now called ``The Selberg zeta-function\footnote{More precisely, the Selberg zeta-function of a cocompact Fuchsian group.}'' as an infinite product over lengths of \emph{primitive closed geodesics}\footnote{Geodesics that do not trace over themselves multiple times.}, bearing a strong resemblance to the Riemann zeta-function. Surprisingly, the Selberg zeta-function satisfies a Riemann hypothosis, and encodes both geometric and spectral data of the quotient orbifold\footnote{A Riemann surface if $\Gamma $ is torsion-free.} $\Gamma \setminus \hh^2.$  The spectral and geometric connection is made clear when one understands the Selberg zeta-function as a by-product of the Selberg trace formula applied to the resolvent kernel of $\lp.$ 

In Defintion~\ref{defSZ} we will define the Selberg zeta-function for our case of interest. A natural question  arises: \emph{what does our zeta-function have in common with  the original Selberg zeta-function?}  The answer\footnote{An alternative answer is that in the cocompact case, both zeta functions are factors of the regularized (functional) determenant $\det \left( \lp - (1-s^2)\right).$  See \cite{Sarnak} and \cite{Friedman} for more details. }: the logarithmic derivatives of both zeta-functions are directly related to the loxodromic (or hyperbolic) contribution of the Selberg trace formula applied to the resolvent kernel of $\lp.$ The term (from the trace formula) in question for our case has the form 
$$
 \sum_{ \{ T \}\LOX}  \frac{   \tr (\chi(T)) \log N(T_{0})}{m(T)|a(T)-a(T)^{-1}|^{2}}N(T)^{-s}.
$$
We show in Lemma~\ref{lemLogDerZ} that the term above is the logarithmic derivatives of a meromorphic function $Z(s,\Gamma,\chi),$ and that it has a product expansion in the right half-plane $R(s)>1.$ 

In order to define $Z(s,\Gamma,\chi)$ we will need some notions concerning centralizer subgroups of loxodromic elements.  For more details see \cite[Sections  5.2,5.4]{Elstrodt}. 

Let  $\Gamma $ be a cofinite  Kleinian group and let $\chi \in \rep.$  Suppose $T\in\Gamma$ is loxodromic (we consider hyperbolic elements as loxodromic elements). Then $T$ is conjugate in $\pc$ to a unique element of the form 
$$
D(T)=
\left(\begin{array}{cc}
a(T) & 0\\
0 & a(T)^{-1}
\end{array}\right) $$
such that $a(T)\in\CC$ has $|a(T)|>1$.  Let $N(T)$  denote  the \emph{norm} of $T,$   defined by  $$N(T) \df |a(T)|^{2},$$ and  let   by  $\mC(T) $ denote  the centralizer of $T$ in $\Gamma.$  There exists a (primitive)  loxodromic element $T_0,$ and a finite cyclic elliptic subgroup  $\en$ of order $m(T), $ generated by an element $E_T $   such that 
$$\mC(T) = \langle T_{0} \rangle \times \en. $$
Here $\langle  T_{0} \rangle = \{\, T_{0}^{n} ~ | ~ n \in\ZZ ~ \}. $ Next,
Let $\mathfrak{t}_1,\dots, \mathfrak{t}_n, $ and $\mathtt{t'_1},\dots,  \mathtt{t'_n}$ denote the eigenvalues of $\chi(T_0)$ and $\chi(E_T)$ respectively.   The elliptic element $ E_T$ is conjugate in $\pc$ to an element of the form 
$$\left(\begin{array}{cc}
\zeta(T_0) & 0 \\
0 & \zeta(T_0)^{-1}
\end{array}\right), $$ 
where here $\zeta(T_0)$ is a primitive $2m(T)$-th root of unity.

\bd \label{defSZ}
For $\R(s)>1 $ the Selberg zeta-function $Z(s,\Gamma,\chi)$ is defined by
$$
Z(s,\Gamma,\chi) \df \prod_{ \{T_0 \} \in \mathcal{R}} ~ \prod_{j=1}^{ \dim V} \prod_{  \substack{ l,k \geq 0 \\  c(T,j,l,k)=1   } } \left( 1-\mathfrak{t}_{j} a(T_0)^{-2k} \overline{ a(T_0) ^{-2l}} N(T_0)^{-s - 1}    \right).  
$$
Here the product with respect to $T_0$ extends over a maximal reduced system $\mathcal{R} $ of $\Gamma$-conjugacy classes of primitive loxodromic elements of $\Gamma.$ The system  $\mathcal{R} $ is called reduced if no two of its elements have representatives with the same centralizer\footnote{See \cite{Elstrodt} section 5.4 for more details}.  The function  $c(T,j,l,k)$ is defined by 
$$c(T,j,l,k)= \mathtt{t'_j} \zeta(T_0)^{2l}  \zeta(T_0)^{-2k}.$$
\ed

\bl \label{lemLogDerZ}For $\R(s)>1,$ 
$$\frac{d}{ds} \log Z(s,\Gamma,\chi)  = \sum_{ \{ T \}\LOX}  \frac{\tr (\chi(T)) \log N(T_{0})}{m(T)|a(T)-a(T)^{-1}|^{2}}N(T)^{-s}. $$
\el
\pf
It follows from the proof of \cite[Lemma 5.4.2]{Elstrodt} that 
\begin{multline}
\sum_{ \{ T \}\LOX}  \frac{\tr (\chi(T)) \log N(T_{0})}{m(T)|a(T)-a(T)^{-1}|^{2}}N(T)^{-s} \\
 = \sum_{ \substack{ \{T_0 \} \in \mathcal{R} \\ n \geq 0 \\ 1 \leq v \leq m(T_0)}}  \frac{ \tr \chi(T_0^{n+1}E_0^v) \log{N(T_0)}  }{m(T_0) \left| \zeta(T_0)^v a(T_0)^{n+1} - \zeta(T_0)^{-v} a(T_0)^{-n-1} \right|^2} N(T_0)^{-s(n+1)}. 
\end{multline}
Next since $T_0$ commutes with $E_T$ we can diagonalize the restriction of $\chi$  to $\mathcal{C}(T) $ and continue the equality to
\begin{multline*}
  = \sum_{ \substack{ \{T_0 \} \in \mathcal{R} \\ n \geq 0 \\ 1 \leq v \leq m(T_0)}} \sum_{j=1}^{\dim V} \frac{ \mathfrak{t}_j^{n+1} \mathtt{t'}_j^v \log{N(T_0)}  }{m(T_0) \left| \zeta(T_0)^v a(T_0)^{n+1} - \zeta(T_0)^{-v} a(T_0)^{-n-1} \right|^2} N(T_0)^{-s(n+1)} \\ 
  = \sum_{ \substack{ \{T_0 \} \in \mathcal{R} \\ n \geq 0 \\ 1 \leq v \leq m(T_0)}} \sum_{j=1}^{\dim V} \frac{ \mathfrak{t}_j^{n+1} \mathtt{t'}_j^v \log{N(T_0)}  }{m(T_0) \left( 1-\zeta(T_0)^{-2v}a(T_0)^{-2(n+1)}  \right) \left( 1-\overline{\zeta(T_0)^{-2v}} \overline{a(T_0)^{-2(n+1)}}  \right) } N(T_0)^{-s(n+1)} \\
  = \sum_{ \substack{ \{T_0 \} \in \mathcal{R} \\ n \geq 0 \\ 1 \leq v \leq m(T_0) \\ l,k \geq 0}} \sum_{j=1}^{\dim V} \frac{ \mathfrak{t}_j^{n+1} \mathtt{t'}_j^v N(T_0)^{-s(n+1)} \log{N(T_0)}  \left( \zeta(T_0)^{-2v}a(T_0)^{-2(n+1)}  \right)^k \left( \overline{\zeta(T_0)^{-2v}} \overline{a(T_0)^{-2(n+1)}}  \right)^l}{m(T_0)}.
\end{multline*}
Next we sum over the $v-$index (note that it is a geometric sum of an $m(T_0)-$th root of unity) observe that the sum is non-zero only when 
$$ \mathtt{t'}_j \zeta(T_0)^{2l}\zeta(T_0)^{-2k}=1 $$ or using our notation $c(T,j,l,k)=1.$   The  equality continues as 
\begin{multline*}
= \sum_{ \substack{ \{T_0 \} \in \mathcal{R} \\ n \geq 0  \\ l,k \geq 0 \\ 1 \leq j \leq \dim V \\ c(T,j,l,k)=1 }}   \mathfrak{t}_j^{n+1}  N(T_0)^{-s(n+1)}\log{N(T_0)}  \left( a(T_0)^{-2(n+1)}  \right)^k \left(  \overline{a(T_0)^{-2(n+1)}}  \right)^l  \\
= \sum_{ \substack{ \{T_0 \} \in \mathcal{R} \\ n,l,k \geq 0 \\ 1 \leq j \leq \dim V \\ c(T,j,l,k)=1}}  \frac{\mathfrak{t}_j a(T_0)^{-2k} \overline{a(T_0)^{-2l}} N(T_0)^{-(s+1)} \log{N(T_0)}}{1-\mathfrak{t}_j a(T_0)^{-2k} \overline{a(T_0)^{-2l}} N(T_0)^{-(s+1)}}
= \frac{Z'(s,\Gamma,\chi)}{Z(s,\Gamma,\chi).}  
\end{multline*}

\epf

\section*{The Logarithmic Derivative of the Selberg Zeta-Function}
The first step in obtaining the meromorphic continuation of the zeta-function is to relate its logarithmic derivative to the trace formula.
From this point on Assumption~\ref{asOne} is in effect. 

We apply the Selberg trace formula to the pair of functions,  $$h(w)=\frac{1}{s^2+w-1} - \frac{1}{B^2+w-1} ~~\text{and}~ $$ 
$$ g(x) = \frac{1}{2s}e^{-s|x|} - \frac{1}{2B}e^{-B|x|}, $$ where $1 < \R(s) < \R(B)$  and obtain

\begin{lem}
\begin{multline} \label{eqLogDer}
\frac{1}{2s} \frac{Z^\prime}{Z}(s) - \frac{1}{2B} \frac{Z^\prime}{Z}(B)  
=\frac{1}{2s} \sum_{ \{ T \}\LOX}  \frac{   \tr (\chi(T)) \log N(T_{0})}{m(T)|a(T)-a(T)^{-1}|^{2}}N(T)^{-s} \\
 -\frac{1}{2B}\sum_{ \{ T \}\LOX} \frac{ \tr (\chi(T))  \log N(T_{0})}{m(T)|a(T)-a(T)^{-1}|^{2}}N(T)^{-B} \\
= \sum_{n \in D} \left(\frac{1}{s^2 - s_n^2} - \frac{1}{B^2 - s_n^2}  \right) 
- \frac{1}{4 \pi} \int_\RR  \left(\frac{1}{s^2 +w^2} - \frac{1}{B^2 + w^2}  \right) \frac{\phi^\prime}{\phi}(i w)~dw  \\
+ \frac{l_\infty}{2 \pi [\gi:\gip]} \int_\RR  \left(\frac{1}{s^2 +w^2} - \frac{1}{B^2 + w^2}  \right) \frac{\Gamma^\prime}{\Gamma}(1 + iw)~dw 
+ \frac{\tr \smat(0)}{4s^2} - \frac{\tr \smat(0)}{4B^2} \\
- \frac{l_\infty}{4[\gi:\gip]s^2} +\frac{l_\infty}{4[\gi:\gip]B^2}  \\
- \sum_{i=1}^l \frac{\tr \chi(g_i) }{|C(g_i)||1-\epsilon_i^2|^2} \int_0^\infty  \left( \frac{e^{-sx}}{2s} - \frac{e^{-Bx}}{2B} \right)   \frac{\sinh x}{\cosh x -1 +\frac{|1-\epsilon_i^2|^2}{2} }~dx 
   \\ - \left( \frac{1}{2s}-\frac{1}{2B}\right) \sum_{ \{R \} \text{\emph{nce}}}\frac{\tr _{V}\chi (R) \log N(T_{0})}{4|\ren |\sin ^{2}(\frac{\pi k}{m(R)})}
+ \frac{\vol \left( \Gamma \setminus \HH \right)\dim V  }{4\pi}(s-B)
 \\ - \left( \frac{1}{2s}-\frac{1}{2B}\right) \sum_{i=1}^l \frac{2 \tr \chi(g_i) \log|c_i| }{|C(g_i)||1-\epsilon_i^2|^2} 
 - \left( \frac{1}{2s}-\frac{1}{2B}\right) \frac{1}{[\gi:\gip]} \left(  l_\infty \left( \frac{ \eta_{\infty}}{2} 
- \gamma \right) +  \sum_{ l = l_\infty + 1 }^{n} L(\Lambda_{\infty}, \psi_{ l}  )   \right).
\end{multline}
\end{lem}
\begin{proof}
The first equality follows from Lemma~\ref{lemLogDerZ}. The second equality follows directly from the Selberg trace formula.
\end{proof}

Equation \eqref{eqLogDer} is  used to exhibit the meromorphic continuation of $Z(s,\Gamma,\chi). $  If we fix $B$ and multiply through   by $2s,$  it is not hard to see that each term on the right of  \eqref{eqLogDer} is meromorphic.  In order to  see that $\sz$ is meromorphic, we must compute the residues of each term on the right of \eqref{eqLogDer}.  We will show that the residues are fractional and that for some $N \in \NN,~\sz^N$ is a meromorphic function.

\begin{thm}
Let  $\Gamma$ be cofinite with one class of cusps at $\zeta = \infty,$ and let $\chi \in \rep.$  
\begin{enumerate}
\item \label{itMer01} If  $[\gi:\gip] = 1$ or $[\gi:\gip] = 2,$  then  $Z(s,\Gamma,\chi)$ is a meromorphic function. 
\item \label{itMer02}If $[\gi:\gip] = 3,$  then there exists a natural number $N,$ $1 \leq N \leq 6,$ so that $\left(Z(s,\Gamma,\chi)\right)^N$ is a meromorphic function.
\end{enumerate}
\end{thm}
\pf
The proof follows from a careful study of \eqref{eqLogDer}.   We must show that after multiplying by $2s,$ each term on the right (of the second equal sign) has at most simple poles with integral or rational residues\footnote{For case (1) the residues must be integer while for case (2) it suffices to show that the residues are rational with bounded denominator. }.   This is demonstrated in  Lemma~\ref{lemSpecDiv}, Lemma~\ref{lemTopCaseOne}, Lemma~\ref{lemTopCaseTwo}, and Lemma~\ref{lemTopCaseThree}.  
\epf
We remark that the divisor of the Selberg zeta-function is readily read off from Lemma~\ref{lemSpecDiv}, Lemma~\ref{lemTopCaseOne}, Lemma~\ref{lemTopCaseTwo}, and Lemma~\ref{lemTopCaseThree}. 

Our zeta function satisfies a functional equation.  A standard argument (\cite[Theorem 5.1.5, page 85]{Venkov}) using \eqref{eqLogDer},  Lemma~\ref{lemSpecDiv}, Lemma~\ref{lemTopCaseOne}, Lemma~\ref{lemTopCaseTwo}, and Lemma~\ref{lemTopCaseThree} yields:
\begin{thm} \label{thmFuncEq}Suppose that $[\gi:\gip]=1$ or $[\gi:\gip]=2.$   Then $\sz$ satisfies:
$$
Z(-s,\Gamma,\chi) = Z(s,\Gamma,\chi)\phi(s)\Psi(s,\Gamma,\chi). 
$$
For $[\gi:\gip]=1,$ 
\begin{equation*}
\Psi(s) \df   \left(\frac{\Gamma(1-s)}{\Gamma(1+s)}   \right)^{k_\infty}  \exp \left(-\frac{\vol \left( \Gamma \setminus \HH \right)\dim V  }{3\pi}s^3+Es + C \right) \end{equation*}
and 
\begin{equation*}
 E \df \sum_{ \{R \} \text{\emph{nce}}}\frac{\tr _{V}\chi (R) \log N(T_{0})}{4|\ren |\sin ^{2}(\frac{\pi k}{m(R)})}
 +  \left(  k_\infty \left( \frac{ \eta_{\infty}}{2} 
- \gamma \right) +  \sum_{ l = k_\infty + 1 }^{n} L(\Lambda_{\infty}, \psi_{ l}  )   \right).
\end{equation*} 
For $[\gi:\gip]=2,$  \begin{multline*}
\Psi(s) \df \\  \left(\frac{\Gamma(1-s)}{\Gamma(1+s)}   \right)^{l_\infty} \left(\prod_{k=1}^\infty \exp(-k(-1)^{k} \frac{(k-1)^2-s^2}{(k+1)^2-s^2}  \right)^{k_\infty/2-l_\infty/2} \exp \left(-\frac{\vol \left( \Gamma \setminus \HH \right)\dim V  }{3\pi}s^3+Es + C \right) \end{multline*}
and 
\begin{multline*}
 E \df \sum_{ \{R \} \text{\emph{nce}}}\frac{\tr _{V}\chi (R) \log N(T_{0})}{4|\ren |\sin ^{2}(\frac{\pi k}{m(R)})}
+  \sum_{i=1}^l \frac{2 \tr \chi(g_i) \log|c_i| }{|C(g_i)||1-\epsilon_i^2|^2} 
 \\ + \frac{1}{[\gi:\gip]} \left(  l_\infty \left( \frac{ \eta_{\infty}}{2} 
- \gamma \right) +  \sum_{ l = l_\infty + 1 }^{n} L(\Lambda_{\infty}, \psi_{ l}  )   \right).
\end{multline*}

 The constant\footnote{The value of $C$ can be read off by letting $s \ra 0$ in the functional equations.  Its value  depends on whether  $\phi(0)$ is $1$ or $-1$ and the multiplicity of $\sz$  at $s=0.$} $C$ satisfies the equation: $\exp(C)=\pm 1.$  
\end{thm}

\bl \label{lemSpecDiv}
The expression 
$$\sum_{n \in D} 2s \left(\frac{1}{s^2 - s_n^2} - \frac{1}{B^2 - s_n^2}  \right) 
- \frac{1}{4 \pi} \int_\RR 2s \left(\frac{1}{s^2 +w^2} - \frac{1}{B^2 + w^2}  \right) \frac{\phi^\prime}{\phi}(i w)~dw   $$ has only simple poles and integral residues:

(a) at the points $\pm s_j$ on the line $\R(s)=0$ and on the interval $[-1,1].$  Each point $s_j$ is related to an eigenvalue $ \lambda_j $ of the discrete spectrum of $\lp $ by $1-s_j^2 = \lambda_j.$  The residue of each $s_j$  is equal to the multiplicity of the corresponding eigenvalue. If $\lambda = 1,$ is an eigenvalue of $\lp, $ then the residue of the point $s_j = 0,$ is twice the multiplicity of $\lambda;$ 

(b) at the points $\rho_j $ that are poles of  $\smat(s),$ which lie in the half-plane $\R(s) < 0. $ The residue of each $\rho_j $ is non-negative\footnote{The point $\rho_j $ is a zero of $Z(s,\Gamma,\chi)$ and a pole of $\smat(s).$ We understand the multiplicity of a pole  as non-negative number (not as a negative number).} and equal to its multiplicity as a pole of $\smat(s).$  
\el
\pf
The computation involves elementary complex analysis. See \cite[Section 5.1]{Venkov}.
\epf

The residues above come  from terms that are related to the spectral and scattering theory of $\lp.$  The remaining residues are computed using group theoretic data involving $\Gamma$ and $\chi.$  The poles and zeros of $\Z(s,\Gamma,\chi)$ that correspond to these residues are commonly called \emph{topological} or \emph{trivial}\footnote{We refer to them as topological.}.

\section*{The Topological Zeros and Poles}
The computation of the topological residues   is considerably more complicated than the corresponding spectral computation.  Poles can only  arise from the following terms (excluding the spectral terms previously dealt with) of \eqref{eqLogDer} (note that we multiplied all terms through by $2s$):
 \begin{multline} \label{eqLogDerTop}
 \frac{l_\infty}{2 \pi [\gi:\gip]} \int_\RR  \frac{2s}{s^2 +w^2}  \frac{\Gamma^\prime}{\Gamma}(1 + iw)~dw  - \frac{l_\infty}{2[\gi:\gip]s}  +  \frac{\tr \smat(0)}{2s}
\\ - \sum_{i=1}^l \frac{\tr \chi(g_i) }{|C(g_i)||1-\epsilon_i^2|^2} \int_0^\infty   e^{-sx}  \frac{\sinh x}{\cosh x -1 +\frac{|1-\epsilon_i^2|^2}{2} }~dx. 
\end{multline}

The first two terms come  from the parabolic elements of $\Gamma$, the third from the spectral trace, and the last from the cuspidal elliptic elements of $\Gamma.$  It is remarkable that the last three terms need to be taken together in order to compute the residue at $s=0,$ while the first and last are needed to compute the residues on the negative real axis.
 
It is well known that 
\beq \label{eqGammaInt}
 \frac{l_\infty}{2 \pi [\gi:\gip]} \int_\RR  \frac{2s}{s^2 +w^2}  \frac{\Gamma^\prime}{\Gamma}(1 + iw)~dw = 
\frac{l_\infty}{[\Gamma_\infty:\Gamma_\infty^\prime]} \left( \frac{\Gamma^\prime}{\Gamma}(1-s) + \sum_{k=1}^\infty \left( \frac{1}{s+k} + \frac{1}{s-k} \right)  \right).
\eeq
In order to obtain a  similarly explicit formula for 
\beq \label{eqCE}
 \sum_{i=1}^l \frac{\tr \chi(g_i) }{|C(g_i)||1-\epsilon_i^2|^2} \int_0^\infty   e^{-sx}  \frac{\sinh x}{\cosh x -1 +\frac{|1-\epsilon_i^2|^2}{2} }~dx
 \eeq
 we must make some technical assumptions.
 \subsection*{Case One: $[\gi:\gip] = 1$} In this case, \eqref{eqCE} is not applicable and  $l_\infty = k_\infty = k(\Gamma,\chi)$ (the last equality follows from our assumption that $\infty$ is the only cusp).  Since $\smat(0)$ is a unitary self-adjoint matrix of dimension $k \times k$, its trace consists of a sum of $k$ terms of the form $ \pm 1.$  It follows that  $\frac{1}{2}(\tr \smat(0) - k)$ is an integer. After applying equations \eqref{eqLogDerTop} and \eqref{eqGammaInt} we have: 
 \bl \label{lemTopCaseOne}
 Suppose $[\gi:\gip] = 1.$ Then the poles of \eqref{eqLogDerTop} are simple and are located at the points $s=-1,-2,\dots,$ with residue $k_\infty$ and at the point $s=0,$ with residue $\frac{1}{2}(\tr \smat(0) - k_\infty).$
 \el
 \subsection*{Case Two: $[\gi:\gip] = 2$ }
 In this case for all $i,~ \epsilon_i = \epsilon \df  \sqrt{-1}  $ and \eqref{eqCE} becomes 
 \beq \label{eqIntTwo}
 \sum_{i=1}^l \frac{\tr \chi(g_i) }{|C(g_i)||1-\epsilon|^2} \int_0^\infty   e^{-sx} \frac{\sinh x}{\cosh x +1  }~dx.
 \eeq
 An application of lemma \ref{lemCuspElip} gives us the coefficient of the integral above 
 $$
 \sum_{i=1}^l \frac{\tr \chi(g_i) }{|C(g_i)||1-\epsilon|^2}  =  \frac{1}{2} \left(
  k_\infty  - \frac{l_\infty }{2}  \right). 
 $$ 
 Next, to evaluate the integral in \eqref{eqIntTwo} we appeal to the following formula:
 \beq
 \int_0^\infty e^{-\mu x} (\cosh x - \cos t)^{-1}~dx = \frac{2}{\sin t} \sum_{k=1}^\infty \frac{\sin kt}{\mu + k},
 \eeq
valid for $\R(\mu) > -1$ and $t \neq 2n\pi, $ (see \cite{Grad} formula  3.543.2).  Averaging for the two values  $\mu -1, \mu+1,$ we obtain 
\beq \label{eqIntegral}
\int_0^\infty e^{-sx}\frac{\sinh x}{\cosh x - \cos t} ~dx  = \frac{1}{\sin t} \sum_{k=1}^\infty \sin{kt} \left( \frac{1}{s-1+k} - \frac{1}{s+1+k} \right).
\eeq
Finally, we can evaluate the integral in \eqref{eqIntTwo} by taking the limit as $t \ra \pi $ 
\beq \label{eqIntTwoEx}
  \int_0^\infty   e^{-sx} \frac{\sinh x}{\cosh x +1  }~dx = 
 \sum_{k=1}^\infty k(-1)^{k+1} \left( \frac{1}{s-1+k} - \frac{1}{s+1+k} \right). 
\eeq
Combining the residues above with those from \eqref{eqGammaInt} give us: 
\bl \label{lemTopCaseTwo}
Suppose $[\gi:\gip]=2.$   Then the poles of \eqref{eqLogDerTop} are simple and are located at the points $s =n,$ with $n = -1,-2,\dots,$ and residue\footnote{Note that by defintion $l_\infty \geq k_\infty.$}:  
$$ m_n = \begin{cases}
k_\infty, 	&  \text{if $n$ is odd}, \\
l_\infty - k_\infty,	& \text{else}
\end{cases}
$$ and at the point $s=0$ with residue $\frac{1}{2}(\tr \smat(0) - k_\infty).$
\el
Since the Picard group $\Gamma = \PSL(2,\ZZ[\sqrt{-1}])$ satisfies $[\gi:\gip]=2$ we have: 
\begin{cor}
Let $\Gamma = \PSL(2,\ZZ[\sqrt{-1}]),$ and let $\chi \in \rep.$ Then $Z(s,\Gamma,\chi)$ is a meromorphic function.
\end{cor}

\subsection*{Case Three: $[\gi:\gip] = 3$ }
In this case, the cuspidal elliptic elements $g_i$ that are in \eqref{eqCE} must all be of order three.  Hence   $\epsilon_i \in \{ \sqrt[3]{1}, \sqrt[3]{-1}~\},$ but when  such $\epsilon_i$ is plugged into $|1-\epsilon_i^2|^2$ one obtains $3.$  Thus we can rewrite \eqref{eqCE} as 
\beq \label{eqCaseThreeA}
\sum_{i=1}^l \frac{\tr \chi(g_i) }{|C(g_i)||1-\epsilon_i^2|^2} \int_0^\infty   e^{-sx}  \frac{\sinh x}{\cosh x + \frac{1}{2} }~dx.
\eeq
Next, applying \eqref{eqIntegral} with $t = \frac{2}{3} \pi, $ and lemma \ref{lemCuspElip}, we can rewrite \eqref{eqCaseThreeA} as
\begin{multline}
\frac{1}{2} \left(k_\infty - \frac{l_\infty}{3} \right) \left( \left( \frac{1}{s-1+1} - \frac{1}{s+1+1} \right) - \left( \frac{1}{s-1+2} - \frac{1}{s+1+2} \right) \right.  \\  \left.
+ \left( \frac{1}{s-1+4} - \frac{1}{s+1+4} \right) - \left( \frac{1}{s-1+5} - \frac{1}{s+1+5} \right) + \dots  \right). 
\end{multline}
Finally  combining the residues above with the  residues from \eqref{eqGammaInt} we obtain: 
\bl \label{lemTopCaseThree}
 Suppose $[\gi:\gip]=3.$   Then the poles of \eqref{eqLogDerTop} are simple and are located at the points $s=n,$ with $n = -1,-2,\dots, $ and residue:
$$ m_n = \begin{cases}
\frac{2}{3}l_\infty - k_\infty, 	&  \text{if $n$ is a multiple of 3}, \\
 \frac{1}{6}l_\infty +\frac{1}{2}k_\infty, 	& \text{else} 
\end{cases}
$$  and at the point $s=0$ with residue $\frac{1}{2}(\tr \smat(0) - k_\infty).$
\el
As an application of the lemma above, we have:
\begin{cor}
Let $\Gamma = \PSL(2,\ZZ[-\frac{1}{2}+\frac{\sqrt{-3}}{2}]),$ and $\chi \equiv 1$ (the trivial representation). Then $Z(s,\Gamma,\chi)$ is not a meromorphic function\footnote{This is the first example that the author is aware of where the Selberg zeta-function is not meromorphic.} (it is the 6-th root of a meromorphic function).
\end{cor}
\pf
Since $\QQ(\sqrt{-3})$ has class number one, $\PSL(2,\ZZ[-\frac{1}{2}+\frac{\sqrt{-3}}{2}])$ has one class of cusps (\cite[Chapter 7]{Elstrodt}).  In addition it has $\infty$ for a cusp.  Elementary calculations show that $[\gi:\gip]=3.$   By definition, $\chi \equiv 1$ implies that $k_\infty = l_\infty =1.$  The result follows from Lemma~\ref{lemTopCaseThree}.
\epf
\subsection*{The Remaining Cases}
The cases that remain are: $[\gi:\gip] = 4$ and $[\gi:\gip] = 6.$  Using the same ideas as in the other three cases one can show the following:
\begin{thm}
Suppose that $[\gi:\gip] = 4.$  Then for some integer $N,$ $\sz^N $ is a meromorphic function.
\end{thm}
The author does not know a good bound for the integer $N$ ($N$ depends on $\Gamma$ and $\chi).$ 
On the other hand we conjecture the following:
\begin{conj*}
Suppose that $[\gi:\gip] = 6.$ Then for some integer $N,$ $\sz^N $ is a meromorphic function.
\end{conj*}
\section*{The Entire Function Associated to the Selberg Zeta-Function}
As Fischer \cite[Chapter 3]{Fischer} observed, it is useful to group the elliptic, parabolic, and identity terms, together with the loxodromic terms to define an entire function associated to the Selberg zeta-function called the \emph{Selberg xi-function} $\Xi(s,\Gamma,\chi).$  
\bd
For $\R(s) > 1 $ set\footnote{Theorem \ref{thmXiFun} shows that we are justified in defining $\Xi(s,\Gamma,\chi)$  by its logarithmic derivative.} 
\begin{multline*}
\frac{\Xi^\prime}{\Xi}(s,\Gamma,\chi) =  \sum_{ \{ T \}\LOX}  \frac{   \tr (\chi(T)) \log N(T_{0})}{m(T)|a(T)-a(T)^{-1}|^{2}}N(T)^{-s} 
+ \frac{l_\infty}{2 \pi [\gi:\gip]} \int_\RR  \left(\frac{2s}{s^2 +w^2}  \right) \frac{\Gamma^\prime}{\Gamma}(1 + iw)~dw \\
+ \frac{1}{2s} \left( \tr \smat(0)
- \frac{l_\infty}{[\gi:\gip]} \right) 
- \sum_{i=1}^l \frac{\tr \chi(g_i) }{|C(g_i)||1-\epsilon_i^2|^2} \int_0^\infty  \left( \frac{e^{-sx}}{2s}  \right)   \frac{\sinh x}{\cosh x -1 +\frac{|1-\epsilon_i^2|^2}{2} }~dx - E,
\end{multline*}
where $E$ is defined in Theorem \ref{thmFuncEq}.
\ed
\begin{thm} \label{thmXiFun}
The Selberg xi-function can be continued to an entire function $\Xi(s,\Gamma,\chi) $ with 
$$
\frac{1}{2s}\frac{\Xi^\prime}{\Xi}(s,\Gamma,\chi) -  \frac{1}{2B}\frac{\Xi^\prime}{\Xi}(B,\Gamma,\chi) = \sum_{n \in D}  \left(\frac{1}{s^2 - s_n^2} - \frac{1}{B^2 - s_n^2}  \right) 
- \frac{1}{4 \pi} \int_\RR  \left(\frac{1}{s^2 +w^2} - \frac{1}{B^2 + w^2}  \right) \frac{\phi^\prime}{\phi}(i w)~dw.   
$$
Here $B > 1$, $B > \R(s).$
\end{thm}
\pf
The second assertion follows from \eqref{eqLogDer}.  For the first asserstion: by Lemma \ref{lemSpecDiv}, after multiplying through by $2s,$ the residues of the expression on the right side of the equal sign are all integers and positive.  Hence $\Xi(s,\Gamma,\chi) $ is entire.
\epf
An explicit product formula (for $\R(s) > 1$) can be obtained for $\Xi(s,\Gamma,\chi) $ in the case of \mbox{$[\gi:\gip]=1,2,3$} by  integrating, and then exponentiating, the explicit formula  for the contribution of the parabolic and cuspidal elliptic elements to the Selberg zeta-function.

\begin{rem} The Selberg zeta-function for symmetric spaces of rank-one was studied by Gangolli and Warner (\cite{Gangolli}, \cite{GangWarn}).  More specifically, they studied torsion-free cocompact quotients with finite-dimensional unitary representations, and non-cocompact torsion-free quotients with  trivial representations (the scalar case).  
\end{rem}
\bibliography{stf}

\providecommand{\bysame}{\leavevmode\hbox to3em{\hrulefill}\thinspace}
\providecommand{\MR}{\relax\ifhmode\unskip\space\fi MR }
\providecommand{\MRhref}[2]{%
  \href{http://www.ams.org/mathscinet-getitem?mr=#1}{#2}
}
\providecommand{\href}[2]{#2}
\begin{thebibliography}{EGM98}

\bibitem[CdV81]{Verdi}
Yves Colin~de Verdi{\`e}re, \emph{Une nouvelle d\'emonstration du prolongement
  m\'eromorphe des s\'eries d'{E}isenstein}, C. R. Acad. Sci. Paris S\'er. I
  Math. \textbf{293} (1981), no.~7, 361--363.

\bibitem[EGM98]{Elstrodt}
J.~Elstrodt, F.~Grunewald, and J.~Mennicke, \emph{Groups acting on hyperbolic
  space}, Springer Monographs in Mathematics, Springer-Verlag, Berlin, 1998,
  Harmonic analysis and number theory.

\bibitem[Fad67]{Faddeev}
L.~D. Faddeev, \emph{The eigenfunction expansion of {L}aplace's operator on the
  fundamental domain of a discrete group on the {L}oba\v cevski\u\i\ plane},
  Trudy Moskov. Mat. Ob\v s\v c. \textbf{17} (1967), 323--350.

\bibitem[Fis87]{Fischer}
J{\"u}rgen Fischer, \emph{An approach to the {S}elberg trace formula via the
  {S}elberg zeta-function}, Lecture Notes in Mathematics, vol. 1253,
  Springer-Verlag, Berlin, 1987.

\bibitem[Fri05]{Friedman}
Joshua Friedman, \emph{The {S}elberg trace formula and {S}elberg zeta-function
  for cofinite {K}leinian groups with finite-dimensional unitary
  representations}, Ph.D. thesis, Stony Brook University, 2005.

\bibitem[Gan77]{Gangolli}
Ramesh Gangolli, \emph{Zeta functions of {S}elberg's type for compact space
  forms of symmetric spaces of rank one}, Illinois J. Math. \textbf{21} (1977),
  no.~1, 1--41. \MR{MR0485702 (58 \#5524)}

\bibitem[GR65]{Grad}
I.~S. Gradshteyn and I.~M. Ryzhik, \emph{Table of integrals, series, and
  products}, Fourth edition prepared by Ju. V. Geronimus and M. Ju. Ce\u\i
  tlin. Translated from the Russian by Scripta Technica, Inc. Translation
  edited by Alan Jeffrey, Academic Press, New York, 1965.

\bibitem[GW80]{GangWarn}
Ramesh Gangolli and Garth Warner, \emph{Zeta functions of {S}elberg's type for
  some noncompact quotients of symmetric spaces of rank one}, Nagoya Math. J.
  \textbf{78} (1980), 1--44.

\bibitem[Hej83]{Hejhal}
Dennis~A. Hejhal, \emph{The {S}elberg trace formula for {${\rm PSL}(2,\,{\bf
  R})$}. {V}ol. 2}, Lecture Notes in Mathematics, vol. 1001, Springer-Verlag,
  Berlin, 1983.

\bibitem[Iwa02]{Iwaniec}
Henryk Iwaniec, \emph{Spectral methods of automorphic forms}, second ed.,
  Graduate Studies in Mathematics, vol.~53, American Mathematical Society,
  Providence, RI, 2002.

\bibitem[Lan76]{Langlands}
Robert~P. Langlands, \emph{On the functional equations satisfied by
  {E}isenstein series}, Springer-Verlag, Berlin, 1976, Lecture Notes in
  Mathematics, Vol. 544.

\bibitem[Lan87]{Lang}
Serge Lang, \emph{Elliptic functions}, second ed., Graduate Texts in
  Mathematics, vol. 112, Springer-Verlag, New York, 1987, With an appendix by
  J. Tate.

\bibitem[Roe66]{Roelcke}
Walter Roelcke, \emph{Das {E}igenwertproblem der automorphen {F}ormen in der
  hyperbolischen {E}bene. {I}, {II}}, Math. Ann. 167 (1966), 292--337; ibid.
  \textbf{168} (1966), 261--324.

\bibitem[Sar87]{Sarnak}
Peter Sarnak, \emph{Determinants of {L}aplacians}, Comm. Math. Phys.
  \textbf{110} (1987), no.~1, 113--120.

\bibitem[Sel56]{Selberg1}
A.~Selberg, \emph{Harmonic analysis and discontinuous groups in weakly
  symmetric {R}iemannian spaces with applications to {D}irichlet series}, J.
  Indian Math. Soc. (N.S.) \textbf{20} (1956), 47--87.

\bibitem[Sel89]{Selberg2}
Atle Selberg, \emph{Collected papers. {V}ol. {I}}, Springer-Verlag, Berlin,
  1989, With a foreword by K. Chandrasekharan.

\bibitem[Sel91]{Selberg3}
\bysame, \emph{Collected papers. {V}ol. {II}}, Springer-Verlag, Berlin, 1991,
  With a foreword by K. Chandrasekharan.

\bibitem[Sie80]{Siegel}
Carl~Ludwig Siegel, \emph{Advanced analytic number theory}, second ed., Tata
  Institute of Fundamental Research Studies in Mathematics, vol.~9, Tata
  Institute of Fundamental Research, Bombay, 1980.

\bibitem[Ven82]{Venkov}
A.~B. Venkov, \emph{Spectral theory of automorphic functions}, Proc. Steklov
  Inst. Math. (1982), no.~4(153), ix+163 pp. (1983), A translation of Trudy
  Mat. Inst. Steklov. {\bf 153} (1981).

\end{thebibliography}
\bibliographystyle{amsalpha}
\end{document}